\documentclass[a4paper, 10pt, twoside, twocolumn]{scrartcl}

\usepackage[utf8]{inputenc} 
\usepackage[T1]{fontenc} 
\usepackage{lmodern}
\usepackage{mathptmx} 
\usepackage[a4paper,twocolumn, twoside]{geometry}         
\usepackage{amsmath} 
\usepackage{amsfonts} 
\usepackage{amssymb} 
\usepackage{graphicx} 
\usepackage{tabularx}
\usepackage[usenames]{color} 
\usepackage{picture}
\usepackage{tikz} 
\usepackage{eso-pic}
\usepackage{fancyhdr} 
\usepackage{calc}
\usepackage{hyperref}
\usepackage{ifthen}
\usepackage{forloop}
\usepackage{setspace}
\usepackage{xstring}
\usepackage[defaultsans]{opensans}
\usepackage{courier}
\newcounter{startpage}

\newcommand{\papertitle}{A new Simheuristics procedure for\\ stochastic combinatorial optimization}

\newcommand{\authors}{Joost Berkhout\textsuperscript{1}}
\newcommand{\institutes}{
     \textsuperscript{1}Department of Mathematics, 
     Vrije Universiteit Amsterdam,\\ {\color{white}\textsuperscript{1}}De Boelelaan 1105, 1081 HV Amsterdam, The Netherlands; {\sffamily\itshape joost.berkhout@vu.nl}
     }
     
\divide\pdfpxdimen by 1000

\definecolor{SNE_DARK_BLUE}{RGB}{23, 54, 93} 
\definecolor{SNE_BLUE}{RGB}{30, 72, 124} 
\definecolor{SNE_BRIGHT_BLUE}{RGB}{54, 94, 144} 
\definecolor{SNE_GRAY}{RGB}{77,77,77} 

\newcounter{ct}
\newcounter{ct2}
\newcommand{\addspaces}[1]{%
  \StrLen{#1}[\abc]
  \setcounter{ct}{\abc}
  \stepcounter{ct}
  \setcounter{ct2}{1}
  \whiledo{\value{ct} > 1}%
  {%
	\StrMid{#1}{\value{ct2}}{\value{ct2}}[\char]%
	\IfEq{\char}{ }{\hspace{1.7mm}}{\char\hspace{1.5mm}}%
	\stepcounter{ct2}%
	\setcounter{ct}{\value{ct}-1}%
   }
}

\newlength{\HSpaceLength}

\makeatletter
\def\hlinewd#1{%
\noalign{\ifnum0=`}\fi\hrule \@height #1 %
\futurelet\reserved@a\@xhline}
\makeatother

\allowdisplaybreaks[4] 

\addtokomafont{subsection}{\color{SNE_BLUE}} 
\addtokomafont{subsection}{\sffamily\bfseries}
\addtokomafont{subsection}{\normalsize}
\addtokomafont{section}{\sffamily\bfseries}
\addtokomafont{section}{\color{SNE_BLUE}}
\addtokomafont{section}{\Large\fontsize{14.1pt}{1em}\selectfont}
\addtokomafont{paragraph}{\sffamily\bfseries}
\addtokomafont{paragraph}{\color{SNE_BLUE}}
\addtokomafont{paragraph}{\normalsize}
\newcommand{\Abstract}[1]{\paragraph{\small\fontsize{8.7pt}{1em}\selectfont\sffamily\bfseries Abstract.}{\small\fontsize{8.7pt}{1em}\selectfont\sffamily #1}}

\addtokomafont{caption}{\sffamily}
\addtokomafont{caption}{\setstretch{1.3}}
\addtokomafont{caption}{\footnotesize\fontsize{8.1pt}{1.2em}\selectfont}
\addtokomafont{caption}{\raggedright}

\renewcommand{\maketitle}{
\twocolumn[\vspace{-3mm}{\sffamily
  {\setstretch{2.36}
    \begin{center}
      {\textcolor{SNE_DARK_BLUE}{\bfseries\huge\fontsize{19.3pt}{1em}\selectfont \papertitle}}\\
      \vspace{-0.4mm}
      {\large\fontsize{12.2pt}{1em}\selectfont \authors}
    \end{center}
  }\vspace{-0.3mm}    
  {\setstretch{1}
    {\small\fontsize{9.1pt}{1em}\selectfont \institutes }
  }}
  \vspace{7.05mm}
]
}
\DeclareMathOperator{\E}{\mathbb{E}}
\newcommand{\Ex}[1]{\E\left[#1\right]}

\definecolor{forestgreen(traditional)}{rgb}{0.0, 0.27, 0.13}

\begin{document}
\maketitle
\noindent
\parbox{\columnwidth}{
\setlength{\fboxrule}{0.8pt}
\setlength{\fboxsep}{1.7mm}
{\color{SNE_GRAY}
\fbox{
  \parbox{0.92\columnwidth}{\setstretch{1.27}\rmfamily\footnotesize\fontsize{8.1pt}{1em}\selectfont
  ~\\
  Here will be SNE-specific data positioned. \\
  To be filled out by editor. \\  
  }
}}
\setstretch{1.1}
\vspace{-1.55 mm}
\Abstract{%
Ignoring uncertainty in combinatorial optimization leads to suboptimal decisions in practice. Nevertheless, the focus is often on deterministic combinatorial optimization problems, mainly because they are already challenging enough without stochasticity. To make it easier to address stochasticity in combinatorial optimization, Simheuristics have been developed that allow solving stochastic combinatorial optimization problems. We propose a new Simheuristic procedure that dynamically changes the optimization focus between a deterministic and stochastic perspective based upon a statistical model. By doing so, an adequate trade-off is made between exploration and exploitation of the solution space during the optimization. We numerically show that the new Simheuristic procedure solves real-life stochastic scheduling problems more efficiently than standard Simheuristics strategies. 
}
}
\normalsize
\setstretch{1.1}
\fontsize{10.3pt}{1.15em}\selectfont

\section{Introduction}

Many real-world problems from various domains, such as logistics, manufacturing, healthcare, and finance, can be stated as combinatorial optimization problems (COPs). These real-life COPs are often NP-hard, meaning there is little hope for an efficient algorithm that allows finding the optimal solution for all realistically-sized problem instances \cite{Papadimitriou1998}. Another complicating factor is that the (input) parameters of COPs are often uncertain in practice \cite{Juan2015}. 

In this work, we focus on these \textit{stochastic} COPs (SCOPs) of the form 
\[
\min_{\pi \in S} \Ex{f(X,\pi)},
\]
where $\pi$ denotes a solution from the discrete solution space $S$, $f(\cdot)$ is the objective function, and $X$ represents the random variable(s) of the parameter(s) with known (empirical) distribution(s). The objective function $f(\cdot)$ is either a closed-form expression or something that can be simulated (for example, a complex production process). A challenging aspect of SCOPs is that $\Ex{f(X,\pi)}$ is generally intractable, and we have to resort to (time-consuming) Monte~Carlo simulations to get sample-average approximations \cite{Juan2015}. In both academia and practice, one often replaces the unknown objective $\Ex{f(X,\pi)}$ by $f(\Ex{X},\pi)$ and optimizes the corresponding \textit{deterministic} COP (DCOP) \cite{Juan2015}:
\[
\min_{\pi \in S} f(\Ex{X},\pi).
\]
Indeed, evaluating a solution $\pi$ in DCOP is done quickly via one function evaluation, whereas finding a good approximation to $\Ex{f(X,\pi)}$ requires many function evaluations. However, this comes at a price that a good solution to DCOP can behave poorly in the corresponding SCOP since $\Ex{f(X,\pi)} \neq f(\Ex{X},\pi)$ in general. Ignoring this is also known as \textit{flaw of averages} \cite{savage_flaw_2012}.

The combination of NP-hardness and the time-consuming objective approximations via simulations make SCOPs challenging to solve in practice. Fortunately, so-called \textit{Simheuristics} have shown to be able to find good solutions to practical SCOPs in recent years \cite{Juan2015}.

\section{Simheuristics}
\begin{figure*}[!ht]
\centering
\mbox{\includegraphics[width=\textwidth]{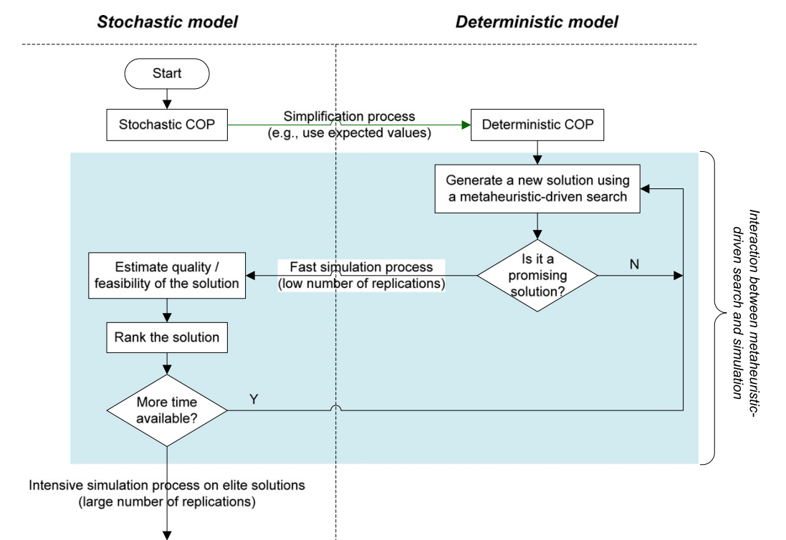}}
\caption{Overview of a Simheuristics framework from \cite{Juan2015}.}
\label{fig:overviewsimheuristics}
\end{figure*}

Simheuristics provide a general framework to solve large-scale SCOPs by combining DCOP (meta)heuristics with simulation \cite{Juan2015}. In particular, the DCOP heuristic is used as a relatively fast way to generate new solutions for SCOP. Instead of simulating all newly found solutions, only the promising solutions are \textit{briefly} simulated to approximate their expected objective values. When approaching the computation time limit, the most promising solutions are awarded more simulations for identifying the best solution to SCOP finally. The process of Simheuristic is illustrated in Figure~\ref{fig:overviewsimheuristics} (taken from \cite{Juan2015}).

Simheuristics are particularly effective for solving SCOPs when: (i) efficient (meta)heuristics already exist for the DCOP, (ii) most gain is obtained in the first part of the DCOP optimization, and (iii) $ f(\Ex{X},\pi)$ and $\Ex{f(X,\pi)}$ are positively correlated for varying $\pi$. As a result, the DCOP heuristic guides the optimization relatively fast to more promising SCOP solutions \cite{Juan2015}. However, when the DCOP optimization stagnates, shifting the optimization focus to simulation is likely more beneficial, i.e., identifying the best SCOP solution out of the most promising DCOP solutions. The key to an effective Simheuristic application is to determine when to switch this optimization focus: switching too early misses out on the chance to find better SCOP solutions efficiently, whereas switching too late increases the chance of picking poor SCOP solutions.

\section{Proposal: OCBA-guided Simheuristic}
We propose a generic procedure for Simheuristics, called \textit{OCBA-guided Simheuristic}, that dynamically determines when to focus on optimizing DCOP and when to focus on simulation to obtain better expected objective values approximations. The idea is to keep track of a fixed-sized elite set of the most promising solutions. At any time during the optimization, we want to be ``sure'' about their expected objective values. To that end, we want the expected opportunity cost of the elite set to be smaller than a user-defined threshold at any time. The opportunity cost of the elite set is the difference between the expected objective values of the solution identified as best and the true best solution, and its expectation is calculated efficiently using Bayesian probability theory \cite{Chen2010}. When the expected opportunity cost exceeds the threshold (meaning we are ``unsure'' about the expected objective values), the solutions from the elite set will be simulated (and we temporarily stop the DCOP optimization). The simulation of the elite set is done efficiently by making use of the Optimal Computing Budget Allocation (OCBA) from \cite{Chen2010}. OCBA prescribes how to efficiently allocate simulation budget among different solutions to minimize the expected opportunity cost. Once the expected opportunity cost drops below the threshold, DCOP is optimized again to find new solutions that may replace solutions from the elite set. This process continues iteratively until the computation budget is spent. Then, the best solution from the elite set is returned.

\section{Preliminary numerical results}
\begin{figure*}[!ht]
\centering
\mbox{\includegraphics[width=\textwidth]{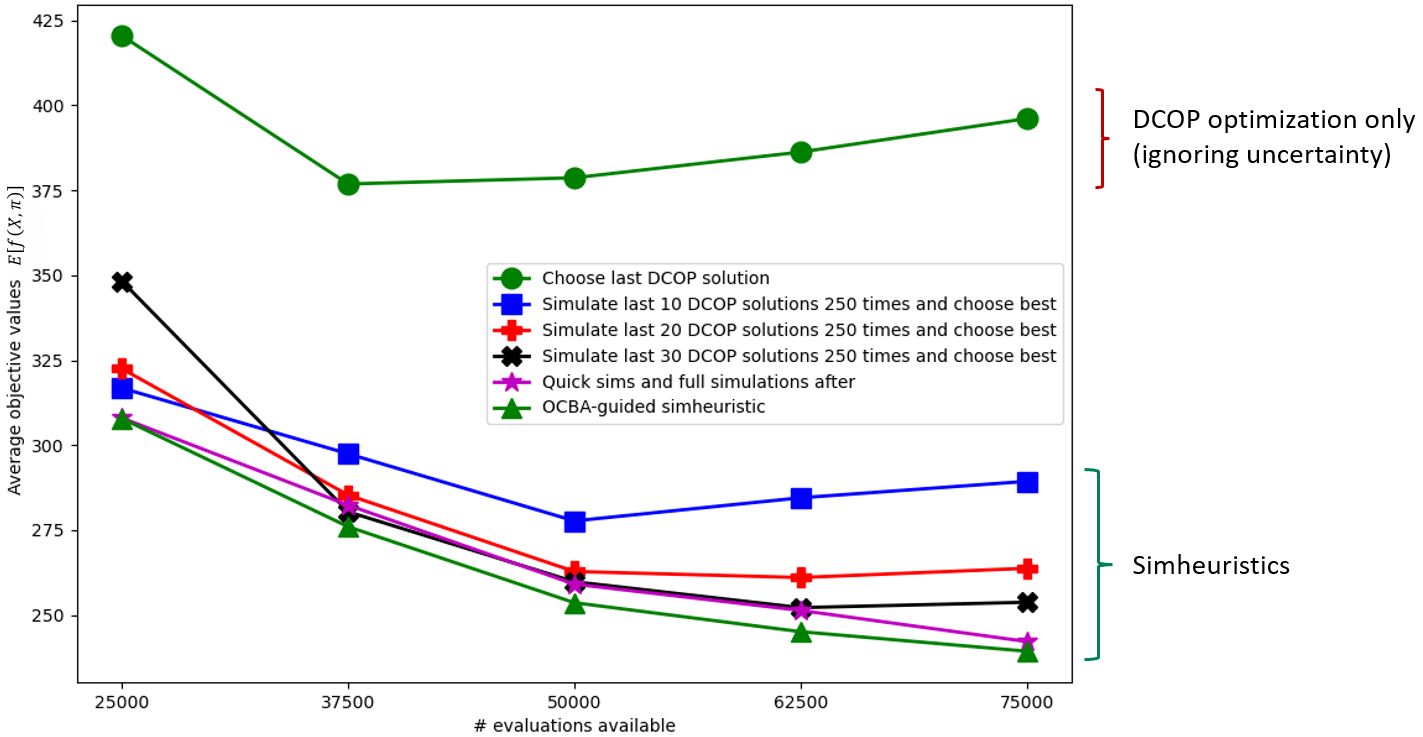}}
\caption{Results of solving a stochastic scheduling problem from \cite{berkhout_short-term_2020} when optimizing DCOP only and using different Simheuristics including our OCBA-guided Simheuristic that dynamically shifts the optimization focus between DCOP and SCOP.}
\label{fig:numerical_results}
\end{figure*}

The OCBA-guided Simheuristic is tested by solving a stochastic version of a parallel machines scheduling problem with sequence-dependent setup times faced in the cattle feed industry \cite{berkhout_short-term_2020}. In particular, we considered for different computation budgets, 50 instances based on real-life data (of 50 jobs, 4 parallel machines, and lognormally distributed production durations) and computed the average expected objective values (which is the weighted sum of the tardiness and the makespan) of the solutions found. The results can be found in Figure~\ref{fig:numerical_results}. For comparison, also the results of optimizing DCOP only and several standard Simheuristics are added. 

\section{Conclusions}
Preliminary results show that our OCBA-guided Simheuristic outperforms other typical Simheuristics for a stochastic scheduling problem. This shows the potential of adequately switching the optimization focus between DCOP and SCOP.

In future research, we want to conduct more experiments. Also, we want to incorporate past simulation information in the OCBA-guided Simheuristic and tailor the simulation budget allocation rule to our purposes.

\vspace{3 mm}

\noindent \textbf{Acknowledgement:} The project has been made possible by TKI Dinalog and the Topsector Logistics and has been funded by the Ministry of Economic Affairs and Climate Policy (EZK).

%
%
%
\bibliographystyle{splncs04}
\bibliography{references}

\end{document}